\documentclass[reqno,10pt]{amsart}
\usepackage{amscd,amssymb,verbatim}
\numberwithin{equation}{section} \theoremstyle{plain}

\newcommand\alp{\alpha}         

\newcommand\gam{\gamma}         \newcommand\Gam{\Gamma}
\newcommand\del{\delta}         \newcommand\Del{\Delta}
\newcommand\eps{\varepsilon}

\newcommand\tet{\theta}

\newcommand\lam{\lambda}

\newcommand\ome{\omega}         \newcommand\Ome{\Omega}

\newcommand\calB{{\mathcal{B}}}
\newcommand\calC{{\mathcal{C}}}

\newcommand\calE{{\mathcal{E}}}

\newcommand\calI{{\mathcal{I}}}

\newcommand\calR{{\mathcal{R}}}

\newcommand\RR{\mathbb{R}}

\newcommand\ZZ{\mathbb{Z}}

\newcommand\CC{\mathbb{C}}

\newcommand\NN{\mathbb{N}}

 \newcommand\gro{{\mathfrak{o}}}

\newcommand\nek{,\ldots,}
\newcommand\sdp{\times \hskip -0.3em {\raise 0.3ex
\hbox{$\scriptscriptstyle |$}}} 

\newcommand\Det{\operatorname{Det}}

\newcommand\Hom{\operatorname {Hom}}
\newcommand\Id{\operatorname {Id}}

\newcommand\Ker{\operatorname{Ker}}

\newcommand\rank{\operatorname{rank}}

\newcommand\Tr{\operatorname{Tr}}

\newcommand\hattau{{\hat{\tau}}}

\renewcommand{\>}{\rangle}
\newcommand{\<}{\langle}

\theoremstyle{plain}
\newtheorem{Thm}[subsection]{Theorem}
\newtheorem{Cor}[subsection]{Corollary}
\newtheorem{Lem}[subsection]{Lemma}
\newtheorem{Prop}[subsection]{Proposition}
\newtheorem{Conjec}[subsection]{Conjecture}

\newtheorem{Def}[subsection]{Definition}

\theoremstyle{remark}

\newtheorem{Rem}[subsection]{Remark}

\def\TeXref#1{%
        \leavevmode\vadjust{\setbox0=\hbox{{\tt
                \  {\tiny \textrm #1}}}%
        \theight=\ht0
        \advance\theight by \lineskip
        \kern -\theight \vbox to
        \theight{\rightline{\rlap{\box0}}%
        \vss}%
        }}%
\newif\ifShowLabels
\ShowLabelstrue
\newdimen\theight
\def\TeXrefEq#1{%
        \leavevmode\vadjust{\setbox0=\hbox{{\tt
                \  {\tiny \textrm #1}}}%
        \theight=\ht1
        \advance\theight by \lineskip
        \kern -\theight \vbox to
        \theight{\rightline{\rlap{\box0}}%
        \vss}%
        }}%

\newcommand{\refs}[1]{Section ~\ref{S:#1}}
\newcommand{\refss}[1]{Subsection ~\ref{SS:#1}}

\newcommand{\reft}[1]{Theorem ~\ref{T:#1}}
\newcommand{\refl}[1]{Lemma ~\ref{L:#1}}
\newcommand{\refp}[1]{Proposition ~\ref{P:#1}}

\newcommand{\refconj}[1]{Conjecture ~\ref{Conj:#1}}
\newcommand{\refd}[1]{Definition ~\ref{D:#1}}

\newcommand{\refe}[1]{\eqref{E:#1}}

\newenvironment{thm}[1]%
        { \begin{Thm} \label{T:#1}  \ifShowLabels \TeXref{T:#1} \fi }%
        { \end{Thm} }

\renewcommand{\th}[1]{\begin{thm}{#1}  }
\renewcommand{\eth}{\end{thm} }

\newenvironment{lemma}[1]%
        { \begin{Lem} \label{L:#1}  \ifShowLabels \TeXref{L:#1} \fi }%
        { \end{Lem} }

\newcommand{\lem}[1]{\begin{lemma}{#1} }
\newcommand{\elem}{\end{lemma}}

\newenvironment{propos}[1]%
        { \begin{Prop} \label{P:#1}  \ifShowLabels \TeXref{P:#1} \fi }%
        { \end{Prop} }

\newcommand{\prop}[1]{\begin{propos}{#1} }
\newcommand{\eprop}{\end{propos}}

\newenvironment{corol}[1]%
        { \begin{Cor} \label{C:#1}  \ifShowLabels \TeXref{C:#1} \fi }%
        { \end{Cor} }
\newcommand{\cor}[1]{\begin{corol}{#1}  }
\newcommand{\ecor}{\end{corol}}

\newenvironment{conjec}[1]%
        { \begin{Conjec} \label{Conj:#1}  \ifShowLabels \TeXref{C:#1} \fi }%
        { \end{Conjec} }
\newcommand{\conj}[1]{\begin{conjec}{#1}  }
\newcommand{\econj}{\end{conjec}}

\newenvironment{defeni}[1]%
        { \begin{Def} \label{D:#1}  \ifShowLabels \TeXref{D:#1} \fi }%
        { \end{Def} }
\newcommand{\defe}[1]{\begin{defeni}{#1}  }
\newcommand{\edefe}{\end{defeni}}

\newenvironment{remark}[1]%
        { \begin{Rem} \label{R:#1}  \ifShowLabels \TeXref{R:#1} \fi }%
        { \end{Rem} }
\newcommand{\rem}[1]{\begin{remark}{#1}}
\newcommand{\erem}{\end{remark}}

\newcommand{\eq}[1]%
        { \ifShowLabels \TeXrefEq{E:#1} \fi
           \begin{equation} \label{E:#1} }
\newcommand{\eeq}{\end{equation}}
\newcommand{\meq}[1]%
        { \ifShowLabels \TeXrefEq{E:#1} \fi
           \begin{multline} \label{E:#1} }
\newcommand{\emeq}{\end{multline}}

\newcommand{\prf}{ \begin{proof} }
\newcommand{\eprf}{ \end{proof} }
\newcommand{\Label}[1]{\label{#1}  \ifShowLabels \TeXref{#1} \fi }

\ShowLabelsfalse

\newcommand{\n}{\nabla}

\renewcommand{\b}{\bullet}
\newcommand{\pa}{\text{\( \partial\)}}

\newcommand{\Arg}{\operatorname{\mathbf{Arg}}}

\newcommand{\Mon}{\operatorname{Mon}}

\newcommand{\Flat}{\operatorname{Flat}}

\newcommand{\B}{\calB}

\newcommand{\rat}{\rho_{\operatorname{an}}}

\newcommand{\trivial}{{\operatorname{trivial}}}

\newcommand{\BH}{{\operatorname{BH}}}
\newcommand{\tauBH}{\tau^\BH}
\newcommand{\Eul}{\operatorname{Eul}}
\newcommand{\PD}{\operatorname{PD}}
\newcommand{\PDp}{\operatorname{PD}'}
\newcommand{\C}{\calC}
\newcommand{\R}{{\operatorname{R}}}

\begin{document}

\title{A Canonical  Quadratic Form on the Determinant Line of a Flat Vector Bundle}
\author[Maxim Braverman]{Maxim Braverman$^\dag$}
\address{Department of Mathematics\\
        Northeastern University   \\
        Boston, MA 02115 \\
        USA
         }
\email{maximbraverman@neu.edu}
\author[Thomas Kappeler]{Thomas Kappeler$^\ddag$}
\address{Institut fur Mathematik\\
         Universitat Z\"urich\\
         Winterthurerstrasse 190\\
         CH-8057 Z\"urich\\
         Switzerland
         }
\email{tk@math.unizh.ch}
\thanks{${}^\dag$Supported in part by the NSF grant DMS-0706837 .\\
\indent${}^\ddag$Supported in part by the Swiss National Science foundation, the programme SPECT, and the European Community through the FP6
Marie Curie RTN ENIGMA (MRTN-CT-2004-5652)}
\begin{abstract}
We introduce and study a canonical quadratic form, called the torsion quadratic form, of the determinant line of a flat vector bundle over a
closed oriented odd-dimensional manifold. This quadratic form caries less information than the refined analytic torsion, introduced in our
previous work, but is easier to construct and closer related to the combinatorial Farber-Turaev torsion. In fact, the torsion quadratic form can
be viewed as an analytic analogue of the Poincar\'e-Reidemeister scalar product, introduced by Farber and Turaev. Moreover, it is also closely
related to the complex analytic torsion defined by Cappell and Miller and we establish the precise relationship between the two. In addition, we
show that up to an explicit factor, which depends on the Euler structure, and a sign the Burghelea-Haller complex analytic torsion, whenever it
is defined, is equal to our quadratic form. We conjecture a formula for the value of the torsion quadratic form at the Farber-Turaev torsion and
prove some weak version of this conjecture. As an application we establish a relationship between the Cappell-Miller and the combinatorial
torsions.
\end{abstract}
\maketitle

\section{Introduction}\Label{S:introduction}

In \cite{BrKappelerBH}, we constructed a new invariant of a flat vector bundle $(E,\n)$ over a closed oriented manifold $M$ of odd dimension $d=2r-1$. It is a
quadratic form $\tau=\tau_\n$, called the {\em torsion quadratic form}, on the determinant line $\Det\big(H^\b(M,E)\big)$ of the cohomology of $E$, which we defined
in terms of another, more sophisticated invariant, the refined analytic torsion $\rat\in \Det\big(H^\b(M,E)\big)$, constructed in
\cite{BrKappelerRAT,BrKappelerRATdetline,BrKappelerRATdetline_hol}.

The invariant $\tau$ is closely related to the quadratic form $\tau^\BH= \tau^\BH_{\n,b}$, introduced by Burghelea and Haller \cite{BurgheleaHaller_function2}. To
construct $\tauBH$ they need to require that the bundle $E$ admits a complex valued non-degenerate bilinear form $b$. The definition of $\tau^{\BH}$ is similar to the
definition of the Ray-Singer torsion, but instead of the standard Laplacians on differential forms uses the non-self-adjoint Laplace-type operators $\Del_b=
\n\n^\#_b+\n^\#_b\n$, where $\n^\#_b$ denotes the adjoint of $\n$ with respect to the bilinear form $b$. Recall that the Ray-Singer torsion is a combination of the
square roots of the determinants of the  standard Laplacians. Since the determinants of the non-self-adjoint operators $\Del_b$ are complex numbers their square roots
are not canonically defined. This is the reason why Burghelea and Haller defined $\tau^\BH$ in terms of the determinants of $\Del_b$ rather than their square roots,
extending in this way the square of the Ray-Singer torsion.

Farber and Turaev \cite{FarberTuraev99,FarberTuraev00} defined a combinatorial torsion $\rho_{\eps,\gro}\in \Det\big(H^\b(M,E)\big)$ which depends on the orientation
$\gro$ of the cohomology $H^\b(M)$ and on the Euler structure $\eps$ introduced by Turaev \cite{Turaev86,Turaev90}. It was noticed by Burghelea \cite{Burghelea99}
that the Euler structure $\eps$ can be described by a closed form $\alp_\eps\in \Ome^{d-1}(M)$. Extending the classical Ray-Singer conjecture,
\cite{RaSi1,Cheeger79,Muller78,BisZh92}, Burghelea and Haller conjectured that
\eq{BHconj}
    \tau^\BH_{\n,b}(\rho_{\eps,\gro}) \ = \ e^{\int_M\ome_{\n,b}\wedge\alp_\eps},
\end{equation}
where $\ome_{\n,b}= -\frac12\Tr{}b^{-1}\n{}b$ is the Kamber-Tondeur form, which measures the non-flatness of the bilinear form $b$. This conjecture was proven
independently by Burghelea-Haller \cite{BurgheleaHaller_function3} and Su-Zhang \cite{SuZhangCM}.

In \cite{BrKappelerBH}, we showed that $\tau^\BH=\pm\tau$ whenever $\tauBH$ is defined and extended the Burghelea-Haller conjecture to the case when $\tauBH$ is not
defined. More precisely, we conjectured, cf. \cite[Conjecture~1.12]{BrKappelerBH}, that
\eq{BKconj}
    \tau_\n(\rho_{\eps,\gro}) \ = \ e^{2\pi i\<\Arg_\n,c(\eps)\>}.
\end{equation}
Here $c(\eps)\in H_1(M,\ZZ)$ is the characteristic class of the Euler structure $\eps$, cf. \cite[\S5.3]{Turaev90}; $\Arg_\n\in H^1(M,\CC/\ZZ)$ is the unique
cohomology class such that for every closed curve $\gam$ in $M$ we have
\eq{IIArgintrod}\notag
    \det\big(\,\Mon_\n(\gam)\,\big) \ = \ \exp\big(\, 2\pi i\<\Arg_\n,[\gam]\>\,\big),
\end{equation}
where $\Mon_\n(\gam)$ denotes the monodromy of the flat connection $\n$ along the curve $\gam$; finally,  $\<\cdot,\cdot\>$ denotes the natural pairing
$H^1(M,\CC/\ZZ)\,\times\, H_1(M,\ZZ) \to \CC/\ZZ$.

Note that \refe{BHconj} implies \refe{BKconj} whenever $\tauBH$ is defined, see \cite[\S1.11]{BrKappelerBH}. In \cite{BrKappelerBH} we proved the following weak
version of Conjecture~\refe{BKconj}: For each connected component $\calC$ of the space of flat connections on $E$ there exists a constant $R_\calC\in \CC$ with
$|R_\calC|=1$, such that
\eq{BKconj2}
    \tau_\n(\rho_{\eps,\gro}) \ = \ R_\calC\cdot e^{2\pi i\<\Arg_\n,c(\eps)\>}.
\end{equation}

Farber and Turaev \cite[\S9]{FarberTuraev00} introduced a bilinear form $\<\cdot,\cdot\>_{\operatorname{PR}}$ on $\Det\big(H^\b(M,E)\big)$, which they call the
cohomological Poincar\'e-Reidemeister scalar product. This is an invariant which refines the Poincar\'e-Reidemeister metric introduced by Farber
\cite{Farber96combinatorialRS}. It follows from Theorem~9.4 of \cite{FarberTuraev00} that Conjecture~\refe{BKconj} is equivalent to the statement that
\[
    \tau_\n(\cdot)\ =\ (-1)^z\, \<\cdot,\cdot\>_{\operatorname{PR}},
\]
where $z\in \NN$ is defined in formula (6.5) of \cite{FarberTuraev00}.

Another related invariant $T\in \Det\big(H^\b(M,E)\big)\otimes \Det\big(H^\b(M,E)\big)$ was introduced by Cappell and Miller \cite{CappellMillerTorsion}. To define
$T$ they also used non-self-adjoint Laplace-type operators, but different from the ones used by Burghelea and Haller. In fact, they consider the square $\B^2$ of the
Atiyah-Patodi-Singer odd signature operator $\B=\B(\n,g^M)$ and, hence, don't need any additional assumptions on $E$. Further  in \cite{CappellMillerTorsion}, Cappell
and Miller conjectured that, in an appropriate sense, their torsion is equal to the Reidemeister torsion of the bundle $E\oplus{}E^*$, where $E^*$ denotes the dual
bundle to $E$.

The goal of this paper is to present a simple construction of the torsion quadratic form $\tau$, implicitly already contained in \cite{BrKappelerRATdetline}. We
collect only those parts of \cite{BrKappelerRAT,BrKappelerRATdetline,BrKappelerBH}, which are needed for this purpose. In particular, we bypass the refined analytic
torsion. Recall that the definition of the refined analytic torsion in \cite{BrKappelerRAT,BrKappelerRATdetline} uses the graded determinant of the odd signature
operator $\B$, leading to a rather complicated analysis, involving the determinant of $\B^2$ and the $\eta$-invariant. In contrast, the definition of $\tau$ only
involves the determinant of the Laplace-type operator $\B^2$. It turns out that the construction of $T$ by Cappell and Miller is very similar to our construction of
$\tau$, as it uses the same operator $\B^2$. We establish the precise relationship of $T$ with $\tau$. It turns out that $T$ is the dual of $\tau$. As an application
we prove a weak version of the Cappell-Miller conjecture.

\section{The Quadratic Form on the Determinant Line of a Finite Dimensional Complex}\label{SS:fdcomplex}

In this section we define a canonical quadratic form on a finite dimensional complex with involution.

\subsection{The construction of a quadratic form}\Label{SS:constrqffd}
Let
\eq{Cpintrod}
    \begin{CD}
       (C^\b,\pa):\quad  0 \ \to C^0 @>{\pa}>> C^1 @>{\pa}>>\cdots @>{\pa}>> C^d \ \to \ 0
    \end{CD}
\end{equation}
be a complex of {\em finite dimensional}\/ complex vector spaces of odd length $d=2r-1$. A {\em chirality operator}\/ $\Gam:C^\b\to C^\b$ is an involution such that
$\Gam(C^j)= C^{d-j}$, for all $j=0\nek d$.  Consider the determinant line
\[
    \Det(C^\b) \ := \ \bigotimes_{j=0}^d\, \Det(C^j)^{(-1)^j},
\]
where $\Det(C^j)^{-1}:= \Hom\big(\Det(C^j),\CC\big)$ denotes the dual of $C^j$. For an element $c_j\in \Det(C^j)$ we denote by $c_j^{-1}$  the unique element in
$\Det(C^j)^{-1}$ satisfying $c_j^{-1}(c_j)=1$. We also denote by $\Gam{}c_j\in \Det(C^{d-j})$ the image of $c_j$ under the map $\Det(C^j)\to \Det(C^{d-j})$ induced by
$\Gam:C^j\to C^{d-j}$.

Denote by $H^\b(\pa)$ the cohomology of the complex $(C^\b,\pa)$. Let
\eq{isomorphism}
    \phi_{C^\b}:\,\Det(C^\b)\ \longrightarrow \  \Det(H^\b(\pa))
\end{equation}
be the canonical isomorphism, cf. \cite{Milnor66}.\footnote{In \cite{BrKappelerRATdetline} we used a sign refined version of this isomorphism,
but we don't need this more complicated version in the present paper.}

Note that any element $c\in \Det(C^\b)$ can be written in a form $c= c_0\otimes c_1^{-1}\otimes \cdots \otimes c_{d}^{-1}$, where $c_j\in \Det(C^j)$. Hence, any
element of $\Det\big(H^\b(\pa)\big)$ can be written as $\phi_{C^\b}(c_0\otimes c_1^{-1}\otimes \cdots \otimes c_{d}^{-1})$.
\defe{taufdim}
The {\em torsion quadratic form} $\tau_{{}_\Gam}$ of the pair $(C^\b,\Gam)$ is the unique quadratic form on $\Det(H^\b(\pa))$ such that
\eq{taufdim}
    \tau_{{}_\Gam}\big(\,\phi_{C^\b}(c_0\otimes c_1^{-1}\otimes \cdots \otimes c_{d}^{-1})\,\big) \ = \ \prod_{j=0}^{d}\,
        \Big[\, c_{j}^{-1}\big(\Gam c_{d-j}\big)\,\Big]^{(-1)^{j+1}}.
\end{equation}
\edefe

\subsection{Relationship with the refined torsion}\label{SS:relreftor}
In \cite{BrKappelerRATdetline} we introduced a canonical element of $\Det\big(H^\b(\pa)\,\big)$, called the {\em refined torsion} of the pair $(C^\b,\Gam)$, as
follows. For each $j=0\nek r-1$, fix an element $c_j\in \Det(C^j)$ and set
\eq{cGam}
    c_{{}_\Gam} \ := \ (-1)^{\calR(C^\b)}\cdot
    c_0\otimes c_1^{-1}\otimes \cdots \otimes c_{r-1}^{(-1)^{r-1}}\otimes (\Gam c_{r-1})^{(-1)^r}\otimes (\Gam c_{r-2})^{(-1)^{r-1}}
    \otimes \cdots\otimes (\Gam c_0)^{-1}
\end{equation}
of $\Det(C^\b)$, where
\eq{R(C)}
  \calR(C^\b) \ = \ \frac12\ \sum_{j=0}^{r-1}\, \dim C^{j}\cdot\big(\, \dim C^j-1\,\big).
\end{equation}
It is easy to see that $c_{{}_\Gam}$ is independent of the choice of $c_0\nek c_{r-1}$. The {\em refined torsion} of the pair $(C^\b,\Gam)$ is the element
\eq{refinedtorintrod}
    \rho_{{}_\Gam} \ = \ \rho_{{}_{C^\b,\Gam}} \ := \  \phi_{C^\b}(c_{{}_\Gam})\ \in \ \Det(H^\b(\pa)).
\end{equation}

It follows immediately from \refe{taufdim} and \refe{refinedtorintrod} that
\eq{tau(rho)}
    \tau_{{}_\Gam}(\rho_{{}_\Gam}) \ = \ 1.
\end{equation}

\subsection{An acyclic complex}\label{SS:tauacyclic}
Suppose the complex $(C^\b,\pa)$ is acyclic. Then $\Det\big(H^\b(\pa)\big)$ is naturally isomorphic to $\CC$. Using this isomorphism we identify $\tau_{{}_\Gam}$ with
the complex number
\eq{tauacyclic}
    \hattau_{{}_\Gam} \ := \ {\tau_{{}_\Gam}(1)} \ \in \CC\backslash\{0\}, \qquad 1\in \CC\simeq \Det\big(H^\b(\pa)\big).
\end{equation}

\subsection{Calculation of the refined torsion of a finite dimensional complex}\Label{SS:calcreftorfd}
To compute the refined torsion we introduce the operator
\eq{Bfd}
    \B \ := \ \Gam\,\pa\ + \pa\,\Gam.
\end{equation}
This operator is a finite dimensional analogue of the signature operator on an odd-dimensional manifold, see \cite[p.~44]{APS1}, \cite[p.~405]{APS2},
\cite[pp.~64--65]{Gilkey84}, and \refs{QuadraticForm} of this paper. Then
\eq{B2}
    \B^2 \ = \ \Gam\,\pa\,\,\Gam\,\pa \ + \ \pa\,\Gam\,\pa\,\Gam.
\end{equation}
\rem{B2=Lapl}
In many interesting applications, cf. \refs{QuadraticForm}, there exists a scalar product on $C^\b$ such that the adjoint of $\pa$ satisfies $\pa^*= \Gam\pa\Gam$.
Then $\B^2$ is equal to the Laplacian of the complex $C^\b$.
\erem

Let us first treat the case where the signature operator $\B$ is bijective.

\lem{tauacyclic}
Suppose that the operator $\B$ is invertible. Then the complex $(C^\b,\pa)$ is acyclic and the complex number $\hattau_{{}_\Gam}$, cf.  \refe{tauacyclic}, is given by
\eq{tauacyclic=}
    \hattau_{{}_\Gam} \ = \ \prod_{j=0}^d\, \Det\big(\,\B^2\big|_{C^j}\,\big)^{(-1)^{j}j}.
\end{equation}
\elem
\prf
Since $\Gam^2=\Id$, for every $a\in \Det(C^\b),\ b\in \Det(C^{d-\b})$, we have
\[
    a^{-1}(\Gam b) \ = \ (\Gam a)^{-1}(b) \ = \ \frac1{b^{-1}(\Gam a)}.
\]
Hence, for all $j=0\nek d$,
\[
    \big[\,c_j^{-1}(\Gam c_{d-j})\,\big]^{(-1)^{j+1}} \ = \ \big[\,c_{d-j}^{-1}(\Gam c_{j})\,\big]^{(-1)^{d-j+1}}
\]
and the definition \refe{taufdim} of $\tau_{{}_\Gam}$ can be rewritten as
\eq{taufdim2}
    \tau_{{}_\Gam}\big(\,\phi_{C^\b}(c_0\otimes c_1^{-1}\otimes \cdots \otimes c_{d}^{-1})\,\big) \ = \ \Big[\,\prod_{j=0}^{r-1}\,
        \Big[\, c_{j}^{-1}\big(\Gam c_{d-j}\big)\,\Big]^{(-1)^{j+1}}\,\Big]^2.
\end{equation}

As, by assumption, the operator $\B=\Gam\pa+\pa\Gam$ is invertible, for each $j=0\nek n$ we have a direct sum decomposition
\[
    C^j \ = \ A^j\oplus B^j,
\]
where $A^j= \Ker\big(\pa|_{C^j}\big)$ and $B^j=\Gam\pa(C^{d-j-1})$. It follows that the complex $(C^\b,\pa)$ is acyclic and $A^j=\pa(B^{j-1})$ for all $j=1\nek d$.
Set $n_j= \dim{}B^j$. Then $n_j= n_{d-j-1}$ and $\dim{}A^j= n_{j-1}$.

For $j=0\nek r-1$ choose a basis $\{b^j_1\nek b^j_{n_j}\}$ of $B^j$. For $j=r\nek d-1$ set $b^j_i= \Gam\pa{}b^{d-j-1}_i$. Then  for any $j=0\nek d-1$, $\{b^j_1\nek
b^j_{n_j}\}$ is a basis of $B^j$. It follows that $\{\pa{}b^{j-1}_1\nek \pa{}b^{j-1}_{n_{j-1}}\}$ is a basis of $A^j$, for $j=1\nek d$. Hence,
\[
    \big\{\,\pa{}b^{j-1}_{1}\nek \pa{}b^{j-1}_{n_{j-1}},b^j_1\nek b^j_{n_j}\,\big\}
\]
is a basis of $C^j$ ($j=1\nek d-1$),\ $\{b^{0}_1\nek b^{0}_{n_{0}}\}$ is the basis of $C^0$, and $\{\pa{}b^{d-1}_1\nek \pa{}b^{d-1}_{n_{d-1}}\}$ is the basis of
$C^d$. Set
\[
    c_0 \ = \ b^{0}_1\wedge\cdots\wedge  b^{0}_{n_{0}}, \qquad c_d \ = \ \pa{}b^{d-1}_1\wedge\cdots\wedge  \pa{}b^{d-1}_{n_{d-1}},
\]
and, for $j=1\nek d-1$,
\[
    c_j \ = \ \pa{}b^{j-1}_{1}\wedge\cdots\wedge \pa{}b^{j-1}_{n_{j-1}}\wedge b^j_1\wedge\cdots\wedge b^j_{n_j} \ \in \ \Det(C^j).
\]
By the definition of the map $\phi_{C^\b}:\Det\big(C^\b\big)\to \Det\big(H^\b(\pa)\big)\simeq \CC$
\[
    \phi_{C^\b}\big(\, c_0\otimes c_1^{-1}\otimes\cdots\otimes c_d^{-1}\,\big) \ = \ 1\ \in \ \CC.
\]
Therefore, by \refe{tauacyclic} and \refe{taufdim2},
\eq{hattau=}
    \hattau_{{}_\Gam} \ = \ \Big[\prod_{j=0}^{r-1}\, \Big[\, c_{j}^{-1}\big(\Gam c_{d-j}\big)\,\Big]^{(-1)^{j+1}}\,\Big]^2.
\end{equation}

We now need to compute the numbers $c_{j}^{-1}\big(\Gam c_{d-j}\big)$.  Assume first, that $j=1\nek r-2$. Then $c_{j}^{-1}\big(\Gam c_{d-j}\big)$ is equal to the
determinant of the operator $S_j:C^j\to C^j$, which transforms the basis $\big\{\pa{}b^{j-1}_{1}\nek \pa{}b^{j-1}_{n_{j-1}},b^j_1\nek b^j_{n_j}\big\}$ to the basis
\eq{newbasis}\notag
    \big\{\,\Gam\pa{}b^{d-j-1}_{1}\nek \Gam\pa{}b^{d-j-1}_{n_{d-j-1}},\Gam b^{d-j}_1\nek \Gam b^{d-j}_{n_{d-j}}\,\big\} \ = \
    \big\{\,\Gam\pa\Gam\pa{}b^{j}_{1}\nek \Gam\pa\Gam\pa{}b^{j}_{n_{j}},\pa b^{j-1}_1\nek \pa b^{j-1}_{n_{j-1}}\,\big\}.
\end{equation}
Here we used that, by construction, $\Gam\pa{}b^j_i= b^{d-j-1}_i$, for any $i=1\nek n_j$ and $b^{d-j}_i= \Gam\pa{}b^{j-1}_i$ for any $i=1\nek n_{j-1}$. We conclude
that
\eq{detSj}
    c_{j}^{-1}\big(\Gam c_{d-j}\big) \ = \ \Det(S_j) \ = \ \pm \Det\big(\,\Gam\pa\Gam\pa\big|_{B^j}\,\big), \qquad j=1\nek r-2.
\end{equation}
Similarly, $c_{0}^{-1}\big(\Gam c_{d}\big)$ is the determinant of the operator which transforms the basis $\{b^0_1\nek b^0_{n_0}\}$ to the basis
\[
    \big\{\,\Gam\pa{}b^{d-1}_{1}\nek \Gam\pa{}b^{d-1}_{n_{d-1}}\,\big\} \ = \
    \big\{\,\Gam\pa\Gam\pa{}b^{0}_{1}\nek \Gam\pa\Gam\pa{}b^{0}_{n_{0}}\,\big\}.
\]
Thus,
\eq{detS0}
    c_{0}^{-1}\big(\Gam c_{d}\big) \ = \  \Det\big(\,\Gam\pa\Gam\pa\big|_{B^0}\,\big).
\end{equation}
Finally, $c_{r-1}^{-1}\big(\Gam c_{r}\big)$ is equal to the determinant of the operator which transforms the basis \linebreak$\big\{\pa{}b^{r-2}_{1}\nek
\pa{}b^{r-2}_{n_{r-2}},b^{r-1}_1\nek b^{r-1}_{n_{r-1}}\big\}$ to the basis
\eq{newbasisr1}\notag
    \big\{\,\Gam\pa{}b^{r-1}_{1}\nek \Gam\pa{}b^{r-1}_{n_{r-1}},\Gam b^{r}_1\nek \Gam b^{r}_{n_{r}}\,\big\} \ = \
    \big\{\,\Gam\pa{}b^{r-1}_{1}\nek \Gam\pa{}b^{r-1}_{n_{r-1}},\pa b^{r-2}_1\nek  \pa b^{r-2}_{n_{r-2}}\,\big\},
\end{equation}
and, hence, is equal to $\pm\Det\big(\Gam\pa\big|_{B^{r-1}}\big)$. Therefore,
\eq{detSr1}
    \big[\,c_{r-1}^{-1}\big(\Gam c_{r}\big)\,\big]^2 \ = \  \Det\big(\,\Gam\pa\big|_{B^{r-1}}\,\big)^2\ = \ \Det\big(\,\Gam\pa\Gam\pa\big|_{B^{r-1}}\,\big).
\end{equation}

Combining equations \refe{hattau=}--\refe{detSr1} we obtain
\eq{hattau=2}
    \hattau_{{}_\Gam} \ = \ \Big[\prod_{j=0}^{r-2}\, \Big[\,\Det\big(\,\Gam\pa\Gam\pa\big|_{B^j}\,\big) \,\Big]^{(-1)^{j+1}}\,\Big]^2\cdot
    \Det\big(\,\Gam\pa\Gam\pa\big|_{B^{r-1}}\,\big).
\end{equation}

The isomorphism $\Gam\pa: B^j\to B^{d-j-1}$ intertwines the operators $\Gam\pa\Gam\pa\big|_{B^j}$ and $\Gam\pa\Gam\pa\big|_{B^{d-j-1}}$. Hence,
\[
    \Det\big(\,\Gam\pa\Gam\pa\big|_{B^j}\,\big) \ = \ \Det\big(\,\Gam\pa\Gam\pa\big|_{B^{d-j-1}}\,\big)
\]
and \refe{hattau=2} can be rewritten as
\eq{hattau=3}
    \hattau_{{}_\Gam} \ = \ \prod_{j=0}^{d-1}\, \Big[\,\Det\big(\,\Gam\pa\Gam\pa\big|_{B^j}\,\big) \,\Big]^{(-1)^{j+1}}.
\end{equation}

The isomorphism $\pa: B^{j-1}\to A^{j}$ intertwines the operators $\Gam\pa\Gam\pa\big|_{B^{j-1}}$ and $\pa\Gam\pa\Gam\big|_{A^{j}}$. Hence,
\[
    \Det\big(\,\Gam\pa\Gam\pa\big|_{B^{j-1}}\,\big) \ = \ \Det\big(\,\pa\Gam\pa\Gam\big|_{A^{j}}\,\big), \qquad j=1\nek d.
\]
Thus, from \refe{B2}, we conclude that
\[
    \Det\big(\, \B^2\big|_{C^0}\,\big) \ = \ \Det\big(\,\Gam\pa\Gam\pa\big|_{B^0}\,\big), \qquad
    \Det\big(\, \B^2\big|_{C^d}\,\big) \ = \ \Det\big(\,\Gam\pa\Gam\pa\big|_{B^{d-1}}\,\big).
\]
and, for $j=1\nek d-1$,
\[
    \Det\big(\, \B^2\big|_{C^j}\,\big) \ = \ \Det\big(\,\Gam\pa\Gam\pa\big|_{B^j}\,\big)\cdot \Det\big(\,\pa\Gam\pa\Gam\big|_{A^{j}}\,\big)
    \ = \ \Det\big(\,\Gam\pa\Gam\pa\big|_{B^j}\,\big)\cdot \Det\big(\,\Gam\pa\Gam\pa\big|_{B^{j-1}}\,\big)
\]
Therefore,
\eq{prod=}
  \begin{aligned}
    \prod_{j=0}^d\, \Det\big(\,\B^2\big|_{C^j}\,\big)^{(-1)^{j}j} \ &= \ \prod_{j=0}^{d-1}\, \Det\big(\,\Det\big(\,\Gam\pa\Gam\pa\big|_{B^j}\,\big)\,\big)^{(-1)^{j}j}
    \cdot \prod_{j=1}^{d}\, \Det\big(\,\Det\big(\,\Gam\pa\Gam\pa\big|_{B^{j-1}}\,\big)\,\big)^{(-1)^{j}j}
    \\ &= \ \prod_{j=0}^{d-1}\, \Det\big(\,\Det\big(\,\Gam\pa\Gam\pa\big|_{B^j}\,\big)\,\big)^{(-1)^{j+1}}.
  \end{aligned}
\end{equation}
Combining \refe{prod=} and \refe{hattau=3} we obtain \refe{tauacyclic=}.
\eprf

To compute the torsion quadratic form in the case $\B$ is {\em not}\/ bijective, note that, for $j=0\nek d$, the operator $\B^2$ maps $C^j$ into itself. For each
$j=0\nek d$ and an arbitrary interval $\calI$, denote by $C^j_\calI\subset C^j$ the linear span of the generalized eigenvectors of the restriction of $\B^2$ to $C^j$,
corresponding to eigenvalues $\lam$ with $|\lam|\in \calI$. Since both operators, $\Gam$ and $\pa$, commute with $\B$ (and, hence, with $\B^2$),
$\Gam(C^j_\calI)\subset C^{d-j}_\calI$ and $\pa(C^j_\calI)\subset C^{j+1}_\calI$. Hence, we obtain a subcomplex $C^\b_\calI$ of $C^\b$ and the restriction
$\Gam_\calI$ of $\Gam$ to $C^\b_\calI$ is a chirality operator for $C^\b_\calI$. We denote by $H^\b_\calI(\pa)$ the cohomology of the complex
$(C^\b_\calI,\pa_\calI)$.

Denote by $\pa_\calI$ and $\B_\calI$ the restrictions of $\pa$ and $\B$ to $C^\b_\calI$. Then $B_\calI= \Gam_\calI\pa_\calI+\pa_\calI\Gam_\calI$ and one easily shows
(cf. Lemma~5.8 of \cite{BrKappelerRATdetline}) that $(C^\b_\calI,\pa_\calI)$ is acyclic if $0\not\in\calI$.

For each $\lam\ge0$, $C^\b= C^\b_{[0,\lam]}\oplus C^\b_{(\lam,\infty)}$ and $H^\b_{(\lam,\infty)}(\pa)= 0$ whereas $H^\b_{[0,\lam]}(\pa)\simeq
H^\b(\pa)$. Hence, there are canonical isomorphisms
\[
        \Phi_\lam:\, \Det(H^\b_{(\lam,\infty)}(\pa))\ \longrightarrow\ \CC, \qquad \Psi_\lam:\,\Det(H^\b_{[0,\lam]}(\pa))\ \longrightarrow \
        \Det(H^\b(\pa)).
\]
In the sequel, we will write $t$ for $\Phi_\lam(t)\in \CC$.

\lem{taugeneral=}
For every $x\in \Det\big(H^\b(\pa)\big)$  and every $\lam\ge0$ we have
\eq{taugeneral=}
    \tau_{{}_\Gam}(x) \ = \ \Big[\,\prod_{j=0}^d\, \Det\big(\,\B^2_{(\lam,\infty)}\big|_{C^j_{(\lam,\infty)}}\,\big)^{(-1)^{j}j}\,\Big] \cdot
    \tau_{{}_{\Gam_{[0,\lam]}}}\big(\,\Psi^{-1}_\lam(x)\,\big).
\end{equation}
In particular, the right hand side of \refe{taugeneral=} is independent of $\lam\ge0$.
\elem
\prf
For each $j=0\nek d$ fix $c_j'\in \Det(C^j_{[0,\lam]})$ and $c_j''\in \Det(C^j_{(\lam,\infty)})$. Then, using the natural isomorphism
\[
    \Det (C^j_{[0,\lam]})\otimes \Det(C^j_{(\lam,\infty)}) \ \simeq \ \Det\big(\,C^j_{[0,\lam]})\oplus C^j_{(\lam,\infty)})\,\big) \ = \ \Det(C^j),
\]
we can regard the tensor product $c_j:=c_j'\otimes{}c_j''$ as an element of $\Det(C^j)$. Applying  \refe{taufdim} twice, we obtain
\meq{tau=tautimestau}
    \tau_{{}_\Gam}\big(\,\phi_{C^\b}(c_0\otimes c_1^{-1}\otimes \cdots \otimes c_{d}^{-1})\,\big)
    \ = \ \prod_{j=0}^{d}\,\Big[\, c_{j}^{-1}\big(\Gam c_{d-j}\big)\,\Big]^{(-1)^{j+1}}
    \\ = \  \prod_{j=0}^{d}\,\Big[\, (c'_{j})^{-1}\big(\Gam c'_{d-j}\big)\,\Big]^{(-1)^{j+1}}\cdot
      \prod_{j=0}^{d}\,
        \Big[\, (c''_{j})^{-1}\big(\Gam c''_{d-j}\big)\,\Big]^{(-1)^{j+1}}
    \\ = \ \tau_{{}_{\Gam_{[0,\lam]}}}\big(\,\phi_{C^\b_{[0,\lam]}}(c'_0\otimes (c'_1)^{-1}\otimes \cdots \otimes (c'_{d})^{-1})\,\big)
        \cdot \tau_{{}_{\Gam_{(\lam,\infty)}}}\big(\,\phi_{C^\b_{(\lam,\infty)}}(c''_0\otimes (c''_1)^{-1}\otimes \cdots \otimes (c''_{d})^{-1})\,\big).
\end{multline}
Let us now choose $c'_j$ and $c''_j$ ($j=0\nek d$) such that $\phi_C^\b(c_0\otimes c_1^{-1}\otimes \cdots \otimes c_{d}^{-1})= x$ and
\[
    \Phi_\lam\circ \phi_{C^\b_{(\lam,\infty)}}\big(\,c''_0\otimes (c''_1)^{-1}\otimes \cdots \otimes (c''_{d})^{-1}\,\big)\ =\ 1.
\]
Then
\[
    \Psi_\lam\circ \phi_{C^\b_{[0,\lam]}}\big(\,c'_0\otimes (c'_1)^{-1}\otimes \cdots \otimes (c'_{d})^{-1}\,\big) \ = \ \pm{}x
\]
and from \refe{tauacyclic=} we get
\[
    \tau_{{}_{\Gam_{(\lam,\infty)}}}\circ\phi_{C^\b_{(\lam,\infty)}}\big(\, c''_0\otimes (c''_1)^{-1}\otimes \cdots \otimes (c''_{d})^{-1}\,\big) \ = \
    \prod_{j=0}^d\, \Det\big(\,\B^2_{(\lam,\infty)}\big|_{C^j_{(\lam,\infty)}}\,\big)^{(-1)^{j}j}.
\]
Hence, \refe{taugeneral=} is equivalent to \refe{tau=tautimestau}.
\eprf

\section{The Quadratic Form Associated to the Square of the Odd Signature Operator}\Label{S:QuadraticForm}

Let $E\to M$ be a complex vector bundle over a closed oriented manifold of {\em odd} dimension $d=2r-1$ and let $\n$ be a flat connection on $E$. Further, let
$\Ome^\b(M,E)$ denote the de Rham complex of $E$-valued differential forms on $M$. For a given Riemannian metric $g^M$ on $M$ denote by
\[
    \Gam \ = \ \Gam(g^M):\, \Ome^\b(M,E) \ \longrightarrow \ \Ome^\b(M,E)
\]
the chirality operator (cf. \cite[\S3]{BeGeVe}), defined in terms of the Hodge $*$-operator by the formula
\eq{Gam}
    \Gam\, \ome \ := \ i^r\,(-1)^{\frac{k(k+1)}2}\,*\,\ome, \qquad \ome\in \Ome^k(M,E).
\end{equation}
The odd signature operator introduced by Atiyah, Patodi, and Singer \cite{APS1,APS2} (see also \cite{Gilkey84}) is the first order elliptic
differential operator $\B:\Ome^\b(M,E)\to \Ome^\b(M,E)$, given by
\[
    \B \ = \ \B(\n,g^M) \ \overset{\text{Def}}{=} \ \Gam\,\n\ +\ \n\,\Gam.
\]
Note that the operator $\B$ is elliptic and its leading symbol is self-adjoint with respect to any Hermitian metric on $E$. Remark also that $\B^2$ maps $\Ome^j(M,E)$
into itself for every $j=0\nek d$. We denote by $(\B^2)_j$ the restriction of $\B^2$ to $\Ome^j(M,E)$.

For an interval $\calI\subset [0,\infty)$ we denote by $\Ome^j_\calI(M,E)$ the image of $\Ome^j(M,E)$ under the spectral projection of $(\B^2)_j$ corresponding to the
eigenvalues whose absolute value lie in $\calI$. The space $\Ome^j_\calI(M,E)$ contains the span of the generalized eigenforms of $(\B^2)_j$ corresponding to
eigenvalues whose absolute value lies in $\calI$ and coincides with this span if the interval $\calI$ is bounded. In particular, since $\B$ is elliptic, if $\calI$ is
bounded, then  the dimension of $\Ome^j_\calI(M,E)$ is finite. Since $\B^2$ and $\n$ commute,  $\Ome^\b_\calI(M,E)$ is a subcomplex of the de Rham complex
$\Ome^\b(M,E)$.

For each $\lam\ge0$, we have
\eq{Ome=Ome0+Ome>0introd}\notag
    \Ome^\b(M,E) \ = \  \Ome^\b_{[0,\lam]}(M,E)\,\oplus\,\Ome^\b_{(\lam,\infty)}(M,E).
\end{equation}
The complex $\Ome^\b_{(\lam,\infty)}(M,E)$ is clearly acyclic. Hence, the cohomology $H^\b_{[0,\lam]}(M,E)$ of the complex
$\Ome^\b_{[0,\lam]}(M,E)$ is naturally isomorphic to the cohomology $H^\b(M,E)$ of $\Ome^\b(M,E)$. Further, as $\Gam$ commutes with $\B^2$, it
preserves the space $\Ome_{[0,\lam]}(M,E)$ and the restriction $\Gam_{\hskip-1pt{}_{[0,\lam]}}$ of $\Gam$ to this space is a chirality operator
on $\Ome^\b_{[0,\lam]}(M,E)$.

Denote by $\B^2_{\calI,j}$ the restrictions of\/ $\B^2$ to $\Ome^j_\calI(M,E)$. Let $\tet\in (0,2\pi)$ be an Agmon angle for $\B_\calI^2$, cf. \cite{ShubinPDObook},
and denote by $\Det_\tet\big(\B^2_{(\lam,\infty),j}\big)$ the $\zeta$-regularized determinant of the operator $\B^2_{(\lam,\infty),j}$ defined using the Agmon angle
$\tet$. Since the leading symbol of $\B^2_{(\lam,\infty),j}$ is positive definite this determinant is independent of the choice of $\tet$.

For any $0\le \lam\le \mu<\infty$, one easily sees that
\eq{DetlamDetmu}
    \prod_{j=0}^d\, \Det_\tet\big(\,\B^2_{(\lam,\infty),j}\,\big)^{(-1)^{j}j} \ = \
    \Big[\,\prod_{j=0}^d\, \Det_\tet\big(\,\B^2_{(\lam,\mu],j}\,\big)^{(-1)^{j}j}\,\Big]\cdot
    \Big[\,\prod_{j=0}^d\, \Det_\tet\big(\,\B^2_{(\mu,\infty),j}\,\big)^{(-1)^{j}j}\,\Big]
\end{equation}

For any given $\lam\ge0$, denote by $\tau_{{}_{\Gam_{\hskip-1pt{}_{[0,\lam]}}}}$ the quadratic form on the determinant line of $H^\b_{[0,\lam]}(M,E)$ associated to
the chirality operator $\Gam_{\hskip-1pt{}_{[0,\lam]}}$, cf. \refd{taufdim}. In view of \refe{taugeneral=} and \refe{DetlamDetmu}, the product
\eq{tau(n)}
    \tau \ = \ \tau(\n) \ := \ \Big[\,\prod_{j=0}^d\, \Det\big(\,\B^2_{(\lam,\infty),j}\,\big)^{(-1)^{j}j}\,\Big]
    \cdot \tau_{{}_{\Gam_{\hskip-1pt{}_{[0,\lam]}}}}
\end{equation}
viewed as a quadratic form on $\Det\big(H^\b(M,E)\big)$ is independent of the choice of $\lam\ge0$. It is also independent of the choice of the Agmon angle $\tet\in
(0,2\pi)$ of $\B^2_{(\lam,\infty)}$.

\defe{taufd}
The quadratic form \refe{tau(n)} on the determinant line of $H^\b(M,E)$ is called the {\em torsion quadratic form}.
\edefe

\th{metric}
The torsion quadratic form $\tau$ is independent of the Riemannian metric $g^M$.
\eth
\prf
Suppose that $g^M_t$, $t\in\RR$, is a smooth family of Riemannian metrics on $M$ and let $\tau_t$ denote the torsion quadratic form corresponding to the metric
$g^M_t$. We need to show that $\tau_t$ is independent of $t$.

Let $\Gam_t$ denote the chirality operator corresponding to the metric $g_t^M$, cf. \refe{Gam}, and let $\B(t)= \B(\n,g^M_t)$ denote the odd signature operator
corresponding to $\Gam_t$.

Fix $t_0\in \RR$ and choose $\lam\ge 0$ so that there are no eigenvalues of $\B(t_0)^2$ whose absolute values are equal to $\lam$. Then there exists $\del>0$ such
that the same is true for all $t\in (t_0-\del,t_0+\del)$. In particular, if we denote by $\Ome^\b_{[0,\lam],t}(M,E)$ the span of the generalized eigenvectors of
$\B(t)^2$ corresponding to eigenvalues with absolute value $\le\lam$, then $\dim\Ome^\b_{[0,\lam],t}(M,E)$ is independent of $t\in (t_0-\del,t_0+\del)$.

Let $\rho_{{}_{\Gam_{\hskip-1pt{}t,[0,\lam]}}}$ denote the refined torsion of the pair $\big(\Ome^\b_{[0,\lam],t}(M,E),\Gam_t\big)$, cf. \refss{relreftor}. As above
we shall view $\rho_{{}_{\Gam_{\hskip-1pt{}t,[0,\lam]}}}$ as an element of $\Det\big(H^\b(M,E)\big)$ via the canonical isomorphism between $H^\b(M,E)$ and
$H^\b_{[0,\lam]}(M,E)$.

In \cite{BrKappelerRATdetline} we fixed a particular square root of $\prod_{j=0}^d\, \Det_\tet\big(\,\B(t)_{(\lam,\infty),j}^2\,\big)^{(-1)^{j+1}{j}}$ (In
\cite{BrKappelerRATdetline} it is denoted by  $e^{\xi_\lam(t,\tet_0)}$. By Lemma~9.2 of \cite{BrKappelerRATdetline} the element
\eq{rho(t)}
    \rho \ := \ \sqrt{\prod_{j=0}^d\, \Det_\tet\big(\,\B(t)_{(\lam,\infty),j}^2\,\big)^{(-1)^{j+1}{j}}}\cdot \rho_{{}_{\Gam_{\hskip-1pt{}t,[0,\lam]}}} \
    \in \ \Det\big(\,H^\b(M,E)\,\big).
\end{equation}
is independent of $t\in (t_0-\del,t_0+\del)$.

Let $\tau_{{}_{\Gam_{\hskip-1pt{}t,[0,\lam]}}}$ denote the torsion quadratic form of the pair  $\big(\Ome^\b_{[0,\lam],t}(M,E),\Gam_t)$. By \refe{tau(rho)} we have
\eq{tau(rho,rho)}
    \tau_t(\rho) \ = \ \prod_{j=0}^d\, \Det\big(\,\B(t)_{(\lam,\infty),j}^2\,\big)^{(-1)^{j}j}\cdot
    \tau_{{}_{\Gam_{\hskip-1pt{}t,[0,\lam]}}}(\rho)
    \ = \ \tau_{{}_{\Gam_{\hskip-1pt{}t,[0,\lam]}}}(\rho_{{}_{\Gam_{\hskip-1pt{}t,[0,\lam]}}}) \ = \ 1,
\end{equation}
where in the latter equality we used \refe{tau(rho)}. Thus $\tau_t(\rho)$ is independent of $t\in (t_0-\del,t_0+\del)$. Since this is true for an arbitrary value of
$t_0$ the theorem is proven.
\eprf
\rem{directproof}
One can easily give a direct proof of \reft{metric}, avoiding any references to \cite{BrKappelerRATdetline}. One only needs to repeat most of
the computations of the proof of Lemma~9.2 of \cite{BrKappelerRATdetline}. However, to save space we preferred to use this lemma, rather than
repeat its proof.
\erem

\section{The Relationship with Burghelea-Haller and Farber-Turaev Torsions}\label{S:BHFT}

In this section we show that the torsion quadratic form $\tau$ coincides with the quadratic form defined in \cite{BrKappelerBH} and use the results of
\cite{BrKappelerBH} to establish the relationship between $\tau$ and the Burghelea-Haller and Farber-Turaev torsions.

\subsection{Relationship with the refined analytic torsion}\label{SS:relrat}
Let $\eta(\n)= \eta(\n,g^M)$ denote the $\eta$-invariant of the restriction of the odd signature operator $\B= \B(\n,g^M)$ to the even forms, see \cite{Gilkey84},
\cite[\S4]{BrKappelerRAT}, \cite[\S6.15]{BrKappelerRATdetline}, or \cite[\S2.9]{BrKappelerBH} for the definition of the $\eta$-invariant of a non-self-adjoint
operator. Let $\eta_\trivial$ be the $\eta$-invariant of trivial line bundle over $M$. Let $\rat= \rat(\n)\in \Det\big(H^\b(M,E)\big)$ denote the refined analytic
torsion of $(E,\n)$, cf. \cite[Definition~9.8]{BrKappelerRATdetline}.

\prop{tau=tau}
\ \ \(\displaystyle \tau_\n\big(\rat(\n)\big) \ = \ e^{-2\pi i\,\big(\,\eta(\n)-\rank E\cdot\eta_\trivial\,\big)}.\)
\eprop
It follows that the torsion quadratic form $\tau$ coincides with the quadratic form defined by equation (1.1) of \cite{BrKappelerBH}.
\prf
Recall that the element $\rho\in \Det\big(\,H^\b(M,E)\,\big)$ is defined in \refe{rho(t)}. From definition of the refined analytic torsion,
\cite[Definition~9.8]{BrKappelerRATdetline}, and formulae (9-5) and (10-21) of \cite{BrKappelerRATdetline} we conclude that
\[
    \rat(\n)\ = \ \pm\,\rho\cdot e^{-\pi i\,\big(\,\eta(\n)-\rank E\cdot\eta_\trivial\,\big)}.
\]
Hence, the statement of the proposition follows immediately from \refe{tau(rho,rho)}.
\eprf

\subsection{Relationship with the Burghelea-Haller torsion}\label{SS:relBH}
Burghelea and Haller \cite{BurgheleaHaller_function,BurgheleaHaller_function2} have introduced a refinement of the square of the Ray-Singer torsion for a closed
manifold of  arbitrary dimension, provided that the complex vector bundle $E$ admits a non-degenerate complex valued symmetric bilinear form $b$. They defined a
complex valued quadratic form
\eq{tauBH}
    \tauBH\ =\ \tauBH_{b,\n}
\end{equation}
on the determinant line $\Det\big(H^\b(M,E)\big)$, which depends holomorphically on the flat connection $\n$ and is closely related to the square of the Ray-Singer
torsion. We refer the reader to \cite{BurgheleaHaller_function,BurgheleaHaller_function2} for the precise definition of the form $\tauBH$ (see also
\cite[\S3]{BrKappelerBH} for a short review).  Using \refp{tau=tau} we now can reformulate Theorem~1.6 of \cite{BrKappelerBH} as follows:
\th{RAT-BH}
Suppose $M$ is a closed oriented manifold of odd dimension $d=2r-1$ and let $E$ be a complex vector bundle over $M$ endowed with a flat connection $\n$. Assume that
there exists a symmetric bilinear form $b$ on $E$ so that the quadratic form \refe{tauBH} on $\Det\big(H^\b(M,E)\big)$ is defined. Then $\tauBH_{b,\n}= \pm\tau_\n$.
\eth
Note that though the Burghelea-Haller form $\tauBH$ is defined only if $E$ admits a non-degenerate bilinear form $b$, the torsion quadratic form $\tau$ exists without
this additional assumption. Therefore, $\tau$ can be viewed as an extension of $\tauBH$ to the case when the bilinear form $b$ does not exist.

\subsection{Relationship with the Farber-Turaev torsion}\label{SS:relFT}
The complex valued combinatorial torsion has been introduced by Turaev \cite{Turaev86,Turaev90,Turaev01} and, in a more general context, by Farber and Turaev
\cite{FarberTuraev99,FarberTuraev00}. The Farber-Turaev torsion depends on the Euler structure $\eps$ and the orientation $\gro$ of the determinant line of the
cohomology $H^\b(M,\RR)$ of $M$. The set of Euler structures $\Eul(M)$, introduced by Turaev, is an affine version of the integer homology $H_1(M,\ZZ)$ of $M$. It has
several equivalent descriptions \cite{Turaev86,Turaev90,Burghelea99,BurgheleaHaller_Euler}. For our purposes, it is convenient to adopt the definition from Section~6
of \cite{Turaev90}, where an Euler structure is defined as an equivalence class of nowhere vanishing vector fields on $M$ -- see \cite[\S5]{Turaev90} for the
description of the equivalence relation. The Farber-Turaev torsion, depending on $\eps$, $\gro$, and $\n$, is an element of the determinant line
$\Det\big(H^\b(M,E)\big)$, which we denote by $\rho_{\eps,\gro}(\n)$.

Suppose $M$ is a closed oriented odd dimensional manifold. Let $\eps\in \Eul(M)$ be an Euler structure on $M$ represented by a non-vanishing vector field $X$,
$\eps=[X]$. Fix a Riemannian metric $g^M$ on $M$ and let $\Psi(g^M)\in \Ome^{d-1}(TM\backslash\{0\})$ denote the Mathai-Quillen form, \cite[\S7]{MathaiQuillen},
\cite[pp.~40-44]{BisZh92}. Set
\[
        \alp_\eps \ = \ \alp_\eps(g^M) \ := \ X^*\Psi(g^M)\ \in \ \Ome^{d-1}(M).
\]
This is a closed differential form, whose cohomology class $[\alp_\eps]\in H^{d-1}(M,\RR)$ is closely related to the integer cohomology class, introduced by  Turaev
\cite[\S5.3]{Turaev90} and called {\em the characteristic class $c(\eps)\in H_1(M,\ZZ)$ associated to an Euler structure $\eps$}. More precisely, let
$\PD:H_1(M,\ZZ)\to H^{d-1}(M,\ZZ)$ denote the Poincar\'e isomorphism. For $h\in H_1(M,\ZZ)$ we denote by $\PDp(h)$ the image of $\PD(h)$ in $H^{d-1}(M,\RR)$. Then
\eq{IPDceps}
    \PD'\big(\,c([X])\,\big) \ = \ -2\,[\alp_\eps] \ = \  -\,2\,[X^*\Psi(g^M)],
\end{equation}

Burghelea and Haller made a conjecture, \cite[Conjecture~5.1]{BurgheleaHaller_function2}, relating the quadratic form $\tauBH_{b,\n}$ and $\rho_{\eps,\gro}(\n)$,
which extends the Bismut-Zhang theorem \cite{BisZh92}. In \cite[Conjecture~1.12]{BrKappelerBH} we extended this conjecture to the case when $E$ does not admit a
non-degenerate symmetric bilinear form. In view of \refp{tau=tau} this conjecture can be reformulated as follows.

Following Farber \cite{Farber00AT}, we denote by $\Arg_\n$ the unique cohomology class $\Arg_\n\in H^1(M,\CC/\ZZ)$ such that for every closed curve $\gam$ in $M$ we
have
\eq{IIArg}\notag
    \det\big(\,\Mon_\n(\gam)\,\big) \ = \ \exp\big(\, 2\pi i\<\Arg_\n,[\gam]\>\,\big),
\end{equation}
where $\Mon_\n(\gam)$ denotes the monodromy of the flat connection $\n$ along the curve $\gam$ and $\<\cdot,\cdot\>$ denotes the natural pairing
$H^1(M,\CC/\ZZ)\,\times\, H_1(M,\ZZ) \to  \CC/\ZZ$.

\conj{BHnew} Assume that $(E,\n)$ is a flat vector bundle over a closed odd dimensional oriented manifold $M$. Then
\eq{BHconjnew}
    \tau_\n\big(\,\rho_{\eps,\gro}(\n)\,\big) \ = \  e^{2\pi i\<\Arg_\n,c(\eps)\>}.
\end{equation}
 \econj

 The original Burghelea-Haller conjecture was proven independently by Burghelea-Haller \cite{BurgheleaHaller_function3} and Su-Zhang \cite{SuZhangCM}. Using this
 result, \reft{RAT-BH}, and formula (1.12) of \cite{BrKappelerBH} we obtain the following theorem, which establishes \refconj{BHnew} up to sign in the case when $E$
 admits a non-degenerate bilinear form:
 \th{BHnewb}
Suppose $M$ is a closed oriented manifold of odd dimension $d=2r-1$ and let $E$ be a complex vector bundle over $M$ endowed with a flat connection $\n$. Assume that
there exists a symmetric bilinear form $b$ on $E$. Then
\eq{BHconjnew1}
    \tau_\n\big(\,\rho_{\eps,\gro}(\n)\,\big) \ = \  \pm\,e^{2\pi i\<\Arg_\n,c(\eps)\>}.
\end{equation}
 \eth
 Also from \refp{tau=tau} and Theorem~1.14 of \cite{BrKappelerBH} we obtain the following
 \th{BHnewweak}
 (i) \ Under the same assumptions as in \refconj{BHnew}, for each connected component $\C$ of the set\/ $\Flat(E)$ of\/ flat connections on $E$ there exists a
 constant $R_\C\in \CC$ with $|R_\C|=1$, such that
\eq{BHconjweak}
    \tau_\n\big(\,\rho_{\eps,\gro}(\n)\,\big) \ = \ R_\C\cdot e^{2\pi i\<\Arg_\n,c(\eps)\>}, \qquad\text{for all}\quad \n\in \C.
\end{equation}

(ii) \ If the connected component $\C$ contains an acyclic Hermitian connection then\/ $R_\C= 1$, i.e.,
\eq{BHconjweak1}
    \tau_\n\big(\,\rho_{\eps,\gro}(\n)\,\big) \ = \  e^{2\pi i\<\Arg_\n,c(\eps)\>}, \qquad\text{for all}\quad \n\in \C.
\end{equation}
\eth
Note that the proof of \reft{BHnewweak} was obtained in \cite{BrKappelerBH} by much softer methods than those used in the proof of the original Burghelea-Haller
conjecture \cite{BurgheleaHaller_function3,SuZhangCM}.

\section{The Cappell-Miller Torsion}\label{S:CappelMiller}

In this section we first recall the definition of Cappell-Miller torsion
\[
    T\ \in\ \Det\big(H^\b(M,E)\big)\otimes\Det\big(H^\b(M,E)\big)
\]
from \cite{CappellMillerTorsion}, then establish its relationship with the torsion form $\tau$, and finally, under some additional assumptions, express  $T$ in terms
of the Farber-Turaev torsion $\rho_{\eps,\gro}$.

\subsection{The Cappell-Miller torsion of a finite dimensional complex}\label{SS:CappellMillerfd}
Let the complex $(C^\b,\pa)$ and the involution $\Gam$ be as in \refss{constrqffd}. Recall that the element $\rho_{{}_\Gam}\in \Det\big(H^\b(\pa)\big)$ was introduced
in \refe{refinedtorintrod}.

In Section~5 of \cite{CappellMillerTorsion} Cappell and Miller introduced a torsion of a class of finite dimensional complexs, which in case of a complex of odd
length $d=2r-1$ and in the presence of the involution $\Gam$ can be described as
\eq{T}
    T = T_{{}_\Gam} \ := \ \rho_{{}_\Gam}\otimes\rho_{{}_\Gam} \in \Det\big(H^\b(\pa)\big)\otimes \Det\big(H^\b(\pa)\big).
\end{equation}

The torsion quadratic form $\tau_{{}_\Gam}$ defined in \refe{taufdim} can be viewed as an element of
\[
    \Det\big(H^\b(\pa)\big)^*\otimes \Det\big(H^\b(\pa)\big)^* \ \simeq \ \Big(\,\Det\big(H^\b(\pa)\big)\otimes \Det\big(H^\b(\pa)\big)\,\Big)^*.
\]
It follows from \refe{tau(rho)} that $\tau_{{}_\Gam}$ is the dual of $T_{{}_\Gam}$, i.e.
\eq{tau=Tfd}
    \tau_{{}_\Gam}(T_{{}_\Gam}) \ = \ 1.
\end{equation}
In particular, if the complex $(C^\b,\pa)$ is acyclic, then $T$ can be viewed as a complex number via the isomorphism $\Det\big(H^\b(\pa)\big)\simeq \CC$, and in
this case $T=1/\tau$. It follows now from \refl{tauacyclic} that if the operator \refe{Bfd} is invertible, then
\eq{T=}
    T_{{}_\Gam} \ = \ \prod_{j=0}^d\, \Det\big(\,\B^2\big|_{C^j}\,\big)^{(-1)^{j+1}j}.
\end{equation}

\rem{CMoriginal}
In \cite{CappellMillerTorsion} the element $T$ is defined in slightly different terms. However, comparing the construction of $\rho_{{}_\Gam}$ with the construction
of Section~5 of \cite{CappellMillerTorsion} one immediately sees that our definition coincides with the one of Cappell-Miller up to sign. To see that the signs agree
one compares \refe{T=} with formula (5.43) of \cite{CappellMillerTorsion}.
\erem

\subsection{The Cappell-Miller torsion of a flat vector bundle}\label{SS:CappelMiller}
Let $E\to M$ be as in \refs{QuadraticForm}. Fix a Riemannian metric $g^M$ on $M$ and let $\Gam$ denote the chirality operator \refe{Gam}. We
shall use the notation introduced in \refs{QuadraticForm}. In particular, for each subset interval $\calI\subset [0,\infty)$ we denote by
$\Ome^j_\calI(M,E)$ the image of $\Ome^j(M,E)$ under the spectral projection of $\B^2\big|_{C^j}$ corresponding to the eigenvalues whose
absolute value lie in $\calI$. Also $\B_{j,\calI}$ denotes the restriction of $\B$ to $\Ome^j_\calI(M,E)$ and $\Gam_\calI$ denotes the
restriction of $\Gam$ to $\Ome^\b_\calE(M,E)$.

Fix $\lam>0$ and let $T_{{}_{\Gam_{[0,\lam]}}}$ be the Cappell-Miller torsion of the complex $\Ome^j_{[0,\lam]}(M,E)$ corresponding to the chirality operator
$\Gam_{[0,\lam]}$. Via the canonical isomorphism $H^\b_{[0,\lam]}(M,E)\simeq H^\b(M,E)$ we can view $T_{{}_{\Gam_{[0,\lam]}}}$ as an element of
$\Det\big(H^\b(\pa)\big)\otimes \Det\big(H^\b(\pa)\big)$.

\defe{CappellMillertorsion}
Let $\tet\in (0,2\pi)$ be an Agmon angle for the operator $\B_{(\lam,\infty)}^2$. The Cappell-Miller torsion $T_\n$ of the flat vector bundle $(E,\n)$ over a closed
oriented odd-dimensional manifold $M$ is the element
\eq{CappellMillertorsion}
     T_\n \ := \ \Big[\,\prod_{j=0}^d\, \Det_\tet\big(\,\B_{(\lam,\infty),j}^2\,\big)^{(-1)^{j+1}j}\,\Big]\cdot T_{{}_{\Gam_{[0,\lam]}}}
     \ \in \ \Det\big(H^\b(\pa)\big)\otimes \Det\big(H^\b(\pa)\big).
\end{equation}
\edefe
It is shown in \cite[Theorem~7.3]{CappellMillerTorsion} and also follows from \reft{tau=T} below that $T_\n$ is independent of the choice of $\lam$.

From \refe{tau=Tfd}, \refe{CappellMillertorsion}, and the definition \refe{tau(n)} of $\tau_\n$ we obtain the following
\th{tau=T}
\(\displaystyle \tau_\n(T_\n) \ = \ 1. \)
\eth
Hence, \refconj{BHnew} can be reformulated in the form
\eq{CappellMillerconj1}
    T_\n \ = \ e^{-2\pi i\<\Arg_\n,c(\eps)\>}\cdot \rho_{\eps,\gro}(\n)\otimes\rho_{\eps,\gro}(\n).
\end{equation}

Let $E^*$ denote the vector bundle dual to $E$. In particular, the fiber $E^*_x$ of $E^*$ at a point $x\in M$ is the dual vector space $E^*_x= \Hom_\CC(E_x,\CC)$. Let
$\n^*$ denote the connection on $E^*$ dual to $\n$. Then the direct sum bundle $E\oplus E^*$ with the connection $\n\oplus\n^*$ is unimodular and its fibers have even
dimension. Hence, cf., for example, Lemmas~3.2 and 3.3 of \cite{FarberTuraev00}, the Reidemeister torsion
\[
    \rho^\R(\n\oplus\n^*)\ \in\ \Det\big(\,H^\b(M,E\oplus{}E^*)\,\big) \simeq \Det\big(\,H^\b(M,E)\,\big)\otimes \Det\big(\,H^\b(M,E^*)\,\big).
\]
is well defined and is equal to the Farber-Turaev torsion $\rho_{\eps,\gro}(\n\oplus{}\n^*)$. In particular, $\rho_{\eps,\gro}(\n\oplus{}\n^*)$ is independent of
$\eps$ and $\gro$.

Farber and Turaev, \cite[p.~219]{FarberTuraev00}, introduced the duality operator
\[
    D:\,\Det\big(\,H^\b(M,E)\,\big)\ \to\  \Det\big(\,H^\b(M,E^*).
\]
Using the definition of the Poincar\'e-Reidemeister scalar product, cf. pages~206 and 219 of \cite{FarberTuraev00} and Theorem~9.4 of
\cite{FarberTuraev00} we obtain
\[
    \rho_{\eps,\gro}(\n)\otimes D\big(\,\rho_{\eps,\gro}(\n)\,\big) \ = \ (-1)^ze^{2\pi i\<\Arg_\n,c(\eps)\>}\cdot\rho^\R(\n\oplus \n^*),
\]
where $z\in \NN$ is defined in formula (6.5) of \cite{FarberTuraev00}. Hence, \refe{CappellMillerconj1} is equivalent to the following
conjecture, originally made by Cappell and Miller \cite{CappellMillerTorsion}:
\begin{Conjec}[\textbf{Cappell-Miller}]\label{Conj:CMconj}
Assume that $(E,\n)$ is a flat vector bundle over a closed odd dimensional oriented manifold $M$. Then the Cappell-Miller torsion is related to the Reidemeister
torsion by the equation
\eq{CMconj}
  (1\otimes D)\,T_\n \ = \ (-1)^z\,\rho^\R(\n\oplus \n^*),
\end{equation}
where $z\in \NN$ is defined in formula (6.5) of \cite{FarberTuraev00}\footnote{The sign factor $(-1)^z$ is missing in
\cite{CappellMillerTorsion} because of a different sign convention.}
\end{Conjec}
Theorems~\ref{T:BHnewb} and  \ref{T:BHnewweak} give a partial solution of this conjecture. In particular, \reft{BHnewweak} says that \refconj{CMconj} holds up to the
factor $R_C$ and holds exactly in the case when $\n$ belongs to a connected component of the space $\Flat(E)$ which contains an acyclic Hermitian connection.
\reft{BHnewb} states that \refconj{CMconj} holds up to sign if $E$ admits a non-degenerate bilinear form $b$.

\providecommand{\bysame}{\leavevmode\hbox to3em{\hrulefill}\thinspace} \providecommand{\MR}{\relax\ifhmode\unskip\space\fi MR }
\providecommand{\MRhref}[2]{%
  \href{http://www.ams.org/mathscinet-getitem?mr=#1}{#2}
} \providecommand{\href}[2]{#2}

\end{document}